\documentstyle[10pt]{article}
\font\tenof=msym10 

\def\Z{\mbox{\tenof Z}}
\def\R{\mbox{\tenof R}}
\def\C{\mbox{\tenof C}}
\def\id{\mbox{\rm id}}
\def\ss{\scriptstyle}

\newtheorem{definition}{Definition}[section]
\newtheorem{proposition}[definition]{Proposition}

\title{Coloured Hopf Algebras\thanks{Presented at the Workshop on Special
Functions and Differential Equations, Chennai, 13--24 January 1997}}
\author{C. Quesne\thanks{Directeur de recherches FNRS; E-mail: cquesne@ulb.ac.be}
\\  
{\small \sl Physique Nucl\'eaire Th\'eorique et Physique Math\'ematique, 
Universit\'e Libre de Bruxelles,} \\ 
{\small \sl Campus de la Plaine CP229,
Boulevard~du Triomphe, B-1050 Brussels, Belgium}}
\date{ }
\begin{document}
\maketitle

\begin{abstract}
Some new algebraic structures related to the coloured Yang-Baxter equation, and
termed coloured Hopf algebras, are reviewed. Coloured quantum universal
enveloping algebras of Lie algebras are defined in this context. An extension to the
coloured graded Yang-Baxter equation and coloured Hopf superalgebras is also
presented. The coloured two-parameter quantum universal enveloping algebra of
$gl(1/1)$ is considered as an example. 
\end{abstract}
%
%
\section{Introduction}     
In recent years, some integrable models with nonadditive-type solutions
$R^{\lambda,\mu} \ne R(\lambda - \mu)$ of the Yang-Baxter equation~(YBE) have
been discovered~\cite{bazhanov}. The corresponding~YBE
\begin{equation}
  R^{\lambda,\mu}_{12} R^{\lambda,\nu}_{13} R^{\mu,\nu}_{23} = R^{\mu,\nu}_{23}
  R^{\lambda,\nu}_{13} R^{\lambda,\mu}_{12}   \label{eq:colYBE}
\end{equation}
is referred to in the literature as the `coloured'~YBE, the nonadditive (in general
multicomponent) spectral parameters $\lambda$, $\mu$,~$\nu$ being considered
as `colour' indices. Constructing solutions of Eq.~(\ref{eq:colYBE}) has been
achieved by using various approaches (see
e.g.~\cite{burdik,ge,kundu,bonatsos}). It should be stressed that this
coloured~YBE is distinct from the so-called `colour YBE'~\cite{mcanally} arising in
another context, as an extension of the graded~YBE to more general gradings than
that determined by $\Z_2$.\par
%
%
Extending the definitions of quantum groups and quantum universal enveloping
algebras (QUEAs)~\cite{majid} by connecting them to coloured $R$-matrices,
instead of ordinary ones, has received some attention in the literature. Kundu and
Basu-Mallick~\cite{kundu} generalized the Faddeev-Reshetikhin-Takhtajan
(FRT) formalism for some quantizations of $U(gl(2))$ and $Gl(2)$. In the context of
knot theory, Ohtsuki~\cite{ohtsuki} introduced some coloured quasitriangular Hopf
algebras, which are characterized by the existence of a coloured universal
$\cal R$-matrix, and he applied his theory to $U_q(sl(2))$. Bonatsos {\it et
al\/}~\cite{bonatsos} independently considered a rather similar, but nevertheless
distinct generalization for some nonlinear deformation of $U(su(2))$. Recently, we
extended the Drinfeld-Jimbo (DJ) formulation of QUEAs of Lie algebras to coloured
ones~\cite{cq}, by elaborating on the results of Bonatsos {\it et~al}.\par
%
%
It is the purpose of the present contribution to review such a generalization. In 
Sec.~2, we define coloured Hopf algebras in a way that generalizes Ohtsuki's
first attempt, and apply the new concepts to QUEAs of Lie algebras. In Sec.~3, we
extend the definitions to the graded case, and consider as an example the
coloured two-parameter QUEA of the Lie superalgebra~$gl(1/1)$.\par
%
%
\section{Coloured Hopf algebras}    
Let $\left({\cal H}_q, +, m_q, \iota_q, \Delta_q, \epsilon_q, S_q;k\right)$ (or in
short ${\cal H}_q$) be a Hopf algebra over some field~$k$ ($= \C$ or \R), depending
upon some parameters $q$~\cite{majid}. Here $m_q: {\cal H}_q \otimes {\cal H}_q
\to {\cal H}_q$, $\iota_q: k \to {\cal H}_q$, $\Delta_q: {\cal H}_q \to {\cal H}_q
\otimes {\cal H}_q$, $\epsilon_q: {\cal H}_q \to k$, and $S_q: {\cal H}_q \to {\cal
H}_q$ denote the multiplication, unit, comultiplication, counit, and antipode maps
respectively. Whenever $q$ runs over some set~$\cal Q$, called
{\it parameter set\/}, we obtain a set of Hopf algebras ${\cal H} = \{\,{\cal H}_q
\mid q \in {\cal Q}\,\}$. We may distinguish between two cases, according to
whether $\cal Q$ contains a single element (fixed-parameter case) or more than one
element (varying-parameter case).\par
%
%
Let us assume~\cite{cq} that there exists a set of one-to-one linear maps ${\cal G}
= \{\,\sigma^{\nu}: {\cal H}_q \to {\cal H}_{q^{\nu}} \mid q, q^{\nu} \in {\cal Q}, \nu
\in {\cal C}\,\}$, defined for any ${\cal H}_q \in {\cal H}$. They are labelled by some
parameters~$\nu$, called {\it colour parameters\/}, taking values in some
set~$\cal C$, called {\it colour set\/}. The latter may be finite, countably infinite,
or uncountably infinite. Two conditions are imposed on the $\sigma^{\nu}$'s:

\begin{itemize}
\item[(i)] Every $\sigma^{\nu}$ is an algebra isomorphism, i.e.,
\begin{equation}
  \sigma^{\nu} \circ m_q = m_{q^{\nu}} \circ \left(\sigma^{\nu} \otimes
  \sigma^{\nu}\right), \qquad \sigma^{\nu} \circ \iota_q = \iota_{q^{\nu}};
  \label{eq:iso} 
\end{equation}
\item[(ii)] $\cal G$ is a group (called {\it colour group\/}) with respect to the
composition of maps, i.e.,
\begin{eqnarray}
  \forall \nu, \nu' \in {\cal C}, \exists\, \nu'' \in {\cal C}: \sigma^{\nu''} & =
           &\sigma^{\nu'} \circ \sigma^{\nu}: {\cal H}_q \to {\cal H}_{q^{\nu''}} = 
          {\cal H}_{q^{\nu,\nu'}},  \label{eq:compo} \\
  \exists\, \nu^0 \in {\cal C}: \sigma^{\nu^0} & = & \id: {\cal H}_q \to 
          {\cal H}_{q^{\nu^0}} = {\cal H}_q, \\
  \forall \nu \in {\cal C}, \exists\, \nu' \in {\cal C}: \sigma^{\nu'} & =
           &\sigma_{\nu} \equiv \left(\sigma^{\nu}\right)^{-1}: {\cal H}_{q^{\nu}} \to
          {\cal H}_{q}   \label{eq:inverse}.
\end{eqnarray}
In Eqs.~(\ref{eq:compo}) and (\ref{eq:inverse}), $\nu''$ and $\nu'$ will be denoted by
$\nu' \circ \nu$ and $\nu^i$, respectively.
\end{itemize}
\par
%
%
$\cal H$, $\cal C$, and $\cal G$ can be combined into
%
\begin{definition}   \label{def-colmaps}
The maps $\Delta^{\lambda,\mu}_{q,\nu}: {\cal H}_{q^{\nu}} \to {\cal
H}_{q^{\lambda}} \otimes {\cal H}_{q^{\mu}}$, $\epsilon_{q,\nu}: {\cal H}_{q^{\nu}}
\to k$, and $S^{\mu}_{q,\nu}: {\cal H}_{q^{\nu}} \to {\cal H}_{q^{\mu}}$, defined by
\begin{equation}
  \Delta^{\lambda,\mu}_{q,\nu} \equiv \left(\sigma^{\lambda} \otimes
  \sigma^{\mu}\right) \circ \Delta_q \circ \sigma_{\nu}, \qquad
  \epsilon_{q,\nu} \equiv \epsilon_q \circ \sigma_{\nu}, \qquad
  S^{\mu}_{q,\nu} \equiv \sigma^{\mu} \circ S_q \circ \sigma_{\nu},
  \label{eq:colmaps} 
\end{equation}
for any $q \in {\cal Q}$, and any $\lambda$, $\mu$, $\nu \in {\cal C}$, are called
coloured comultiplication, counit, and antipode respectively. 
\end{definition}
\par
%
%
It is easy to prove the following proposition:
%
\begin{proposition}   \label{prop-GenHopf}
The coloured comultiplication, counit, and antipode maps, defined in
Eq.~(\ref{eq:colmaps}), transform under the colour group~$\cal G$ as
\begin{eqnarray}
  \left(\sigma^{\lambda}_{\alpha} \otimes \sigma^{\mu}_{\beta}\right) \circ 
           \Delta^{\alpha,\beta}_{q,\nu} & = & \Delta^{\lambda,\mu}_{q,\nu} = 
           \Delta^{\lambda,\mu}_{q,\gamma} \circ \sigma^{\gamma}_{\nu}, 
           \nonumber \\
  \epsilon_{q,\alpha} \circ \sigma^{\alpha}_{\nu} & = & \epsilon_{q,\nu},
           \nonumber \\
  \sigma^{\mu}_{\alpha} \circ S^{\alpha}_{q,\nu} & = & S^{\mu}_{q,\nu} = 
           S^{\mu}_{q,\beta} \circ \sigma^{\beta}_{\nu},
\end{eqnarray}
and satisfy generalized coassociativity, counit, and antipode axioms
\begin{eqnarray}
  \left(\Delta^{\alpha,\beta}_{q,\lambda} \otimes \sigma^{\gamma}_{\mu}\right)
           \circ \Delta^{\lambda,\mu}_{q,\nu} & = & \left(\sigma^{\alpha}_{\lambda'}
           \otimes \Delta^{\beta,\gamma}_{q,\mu'}\right) \circ 
           \Delta^{\lambda',\mu'}_{q,\nu}, \nonumber \\
  \left(\epsilon_{q,\lambda} \otimes \sigma^{\alpha}_{\mu}\right) \circ
           \Delta^{\lambda,\mu}_{q,\nu} & = & \left(\sigma^{\alpha}_{\lambda'}
           \otimes \epsilon_{q,\mu'}\right) \circ \Delta^{\lambda',\mu'}_{q,\nu} =
           \sigma^{\alpha}_{\nu}, \nonumber \\
  m_{q^{\alpha}} \circ \left(S^{\alpha}_{q,\lambda} \otimes \sigma^{\alpha}_{\mu}
           \right) \circ \Delta^{\lambda,\mu}_{q,\nu} & = & m_{q^{\alpha}} \circ \left(
           \sigma^{\alpha}_{\lambda'} \otimes S^{\alpha}_{q,\mu'} \right) \circ
           \Delta^{\lambda',\mu'}_{q,\nu} = \iota_{q^{\alpha}} \circ \epsilon_{q,\nu},
           \label{eq:colantipode}
\end{eqnarray}
as well as generalized bialgebra axioms
\begin{eqnarray}
  \Delta^{\lambda,\mu}_{q,\nu} \circ m_{q^{\nu}} & = & \left(m_{q^{\lambda}}
           \otimes m_{q^{\mu}}\right) \circ (\id \otimes \tau \otimes \id) \circ
           \left(\Delta^{\lambda,\mu}_{q,\nu} \otimes \Delta^{\lambda,\mu}_{q,\nu}
           \right), \nonumber \\
  \Delta^{\lambda,\mu}_{q,\nu} \circ \iota_{q^{\nu}} & = & \iota_{q^{\lambda}}
           \otimes \iota_{q^{\mu}}, \nonumber \\
  \epsilon_{q,\nu} \circ m_{q^{\nu}} & = & \epsilon_{q,\nu} \otimes \epsilon_{q,\nu},
           \nonumber \\
  \epsilon_{q,\nu} \circ \iota_{q^{\nu}} & = & 1_k.
\end{eqnarray}
Here $\sigma^{\lambda}_{\mu}$ is the element of~$\cal G$ defined by
\begin{equation}
  \sigma^{\lambda}_{\mu} \equiv \sigma^{\lambda} \circ \sigma_{\mu},
\end{equation}
$\tau$ is the twist map, i.e., $\tau(a \otimes b) = b \otimes a$, $1_k$ denotes the
unit of~$k$, and no summation is implied over repeated indices.
\end{proposition}
\par
%
%
We are then led to introduce
%
\begin{definition}  \label{def-colHopf}
A set of Hopf algebras~$\cal H$, endowed with coloured comultiplication, counit,
and antipode maps $\Delta^{\lambda,\mu}_{q,\nu}$, $\epsilon_{q,\nu}$,
$S^{\mu}_{q,\nu}$, as defined in (\ref{eq:colmaps}), is called coloured Hopf algebra,
and denoted by any one of the symbols $\left({\cal H}_q, +, m_q, \iota_q,
\Delta^{\lambda,\mu}_{q,\nu}, \epsilon_{q,\nu}, S^{\mu}_{q,\nu}; k, {\cal Q}, {\cal
C}, {\cal G}\right)$, $\left({\cal H}, {\cal C}, {\cal G}\right)$, or ${\cal H}^c$.
\end{definition}
%
\par
%
%
Let us now assume that the members of the Hopf algebra set~$\cal H$ are
quasitriangular Hopf algebras $\left({\cal H}_q, {\cal R}_q\right)$,
where ${\cal R}_q \in {\cal H}_q \otimes {\cal H}_q$ denotes the corresponding
universal $\cal R$-matrix~\cite{majid}. We may then introduce
%
\begin{definition}    \label{def-colR}
Let ${\cal R}^c$ denote the set of elements ${\cal R}^{\lambda,\mu}_q \in {\cal
H}_{q^{\lambda}} \otimes {\cal H}_{q^{\mu}}$, defined by
\begin{equation}
  {\cal R}^{\lambda,\mu}_q \equiv \left(\sigma^{\lambda} \otimes
  \sigma^{\mu}\right) \left({\cal R}_q\right),    \label{eq:colR}  
\end{equation}
where $q$ runs over~$\cal Q$, and $\lambda$, $\mu$ over~$\cal C$.
\end{definition}
\par
%
%
The following result can be easily obtained:
%
\begin{proposition}
If the Hopf algebras~${\cal H}_q$ of~$\cal H$ are quasitriangular, then
${\cal R}^{\lambda,\mu}_q$, as defined in (\ref{eq:colR}), is invertible with
$\left({\cal R}^{\lambda,\mu}_q\right)^{-1}$ given by
\begin{equation}
  \left({\cal R}^{\lambda,\mu}_q\right)^{-1} = \left(\sigma^{\lambda} \otimes
  \sigma^{\mu}\right) \left({\cal R}_q^{-1}\right),   \label{eq:quasi1}
\end{equation}
and for any $a \in {\cal H}_{q^{\nu}}$,
\begin{eqnarray}
  \tau \circ \Delta^{\mu,\lambda}_{q,\nu}(a) & = & {\cal R}^{\lambda,\mu}_q
            \Delta^{\lambda,\mu}_{q,\nu}(a) \left({\cal R}^{\lambda,\mu}_q\right)^{-1},
            \label{eq:quasi2} \\  
  \left(\Delta^{\alpha,\beta}_{q,\lambda} \otimes \sigma^{\gamma}_{\mu}\right)
            \left({\cal R}^{\lambda,\mu}_q\right) & = & {\cal R}^{\alpha,\gamma}_{q,13}
            {\cal R}^{\beta,\gamma}_{q,23}, \label{eq:quasi3} \\
  \left(\sigma^{\alpha}_{\lambda} \otimes \Delta^{\beta,\gamma}_{q,\mu}\right)
            \left({\cal R}^{\lambda,\mu}_q\right) & = & {\cal R}^{\alpha,\gamma}_{q,13}
            {\cal R}^{\alpha,\beta}_{q,12}.    \label{eq:quasi4}
\end{eqnarray}
\end{proposition}
%
Hence we have
%
\begin{definition}  \label{def-quasi}
A coloured quasitriangular Hopf algebra is a pair $\left({\cal H}^c, {\cal
R}^c\right)$, where ${\cal H}^c$ is a coloured Hopf algebra, ${\cal R}^c = \{\, {\cal
R}^{\lambda,\mu}_q \mid q \in {\cal Q}, \lambda, \mu \in {\cal C}\,\}$, and ${\cal
R}^{\lambda,\mu}_q$, defined in (\ref{eq:colR}), satisfies Eqs.~(\ref{eq:quasi1}),
(\ref{eq:quasi2}), (\ref{eq:quasi3}), and~(\ref{eq:quasi4}). The set~${\cal R}^c$ is 
called the coloured universal $\cal R$-matrix of $\left({\cal H}^c, {\cal R}^c\right)$.
\end{definition}
\par
%
%
The terminology used for~${\cal R}^c$ in Definition~\ref{def-quasi} is justified by
the following proposition:
%
\begin{proposition}
Let $\left({\cal H}^c, {\cal R}^c\right)$ be a coloured quasitriangular Hopf algebra.
Then the elements of~${\cal R}^c$ satisfy the coloured YBE
\begin{equation}
  {\cal R}^{\lambda,\mu}_{q,12} {\cal R}^{\lambda,\nu}_{q,13} 
            {\cal R}^{\mu,\nu}_{q,23} = {\cal R}^{\mu,\nu}_{q,23}
            {\cal R}^{\lambda,\nu}_{q,13} {\cal R}^{\lambda,\mu}_{q,12}.
\end{equation}
\end{proposition}
%
\par
%
%
As shown elsewhere~\cite{cq}, the coloured (quasitriangular) Hopf algebras
introduced in Definition~\ref{def-colHopf} (and Definition~\ref{def-quasi}) 
generalize those previously introduced by Ohtsuki~\cite{ohtsuki}, which are
restricted to abelian colour groups, in which case they reduce to substructures of
the present ones.\par
%
%
Extending the DJ~formulation of QUEAs of Lie algebras, $U_q(g)$, to coloured ones
is now straightforward:
%
\begin{definition}   \label{def-colQUEA}
If ${\cal H}_q = U_q(g)$, where $g$ is some Lie algebra, then the corresponding
coloured Hopf algebra~${\cal H}^c$ is called a coloured QUEA.
\end{definition}
\par
%
%
In Ref.~\cite{cq}, various examples of coloured QUEAs of Lie algebras have been
constructed. In the fixed-parameter case, they correspond to

\begin{itemize}
\item[(i)] the standard quantum algebra~$U_q(sl(2))$, 
\item[(ii)] the two-parameter quantum algebra~$U_{q,s}(gl(2))$, 
\item[(iii)] the three-parameter quantum algebra~$U_{q,s_1,s_2}(sl(3) \oplus u(1)
\oplus u(1))$,
\end{itemize}

\noindent while in the varying-parameter case, they are related to

\begin{itemize}
\item[(i)] the nonstandard quantum algebra~$U_h(sl(2))$, 
\item[(ii)] the standard quantum oscillator algebra~$U^{(s)}_z(h(4))$, 
\item[(iii)] the one-parameter nonstandard quantum oscillator
algebra~$U^{(n)}_z(h(4))$, 
\item[(iv)] the three-parameter nonstandard quantum oscillator
algebra~$U^{(IIn)}_{\vartheta,\beta_+,\beta_-}(h(4))$, 
\item[(v)] the standard three-dimensional quantum Euclidean algebra~$U_w(e(3))$,
\item[(vi)] the null-plane $D$-dimensional quantum Poincar\'e
algebras~$U_z(iso(D-1,1))$. 
\end{itemize}

Such examples show that any QUEA can be easily transformed into a coloured
one, that this can be achieved in various ways, and that some of them may involve a
nonabelian colour group.\par
%
%
\section{Extension to the graded case}
Let us now assume that the Hopf algebras~${\cal H}_q$, considered in the previous
section, are \Z$_2$-graded, and denote by $\gamma_q: {\cal H}_q \to {\cal H}_q$
their grading automorphism, i.e.,
\begin{equation}
  \gamma_q(a) = (-1)^{\deg a} a,
\end{equation}
for any homogeneous $a \in {\cal H}_q$, where $\deg a = 0$ or~1
according to whether $a$ is even or odd. As the spaces are graded, the tensor
product and the twist map are now such that
\begin{equation}
  (a \otimes b) (c \otimes d) = (-1)^{(\deg b)(\deg c)} ac \otimes bd, \qquad
  \tau(a \otimes b) = (-1)^{(\deg a)(\deg a)} b \otimes a,   \label{eq:permute} 
\end{equation}
for any homogeneous $a$, $b$, $c$, $d \in {\cal H}_q$~\cite{chaichian}.\par
%
%
Provided we take these properties into account, it is straightforward to extend the
definitions of the previous section in a way that preserves the grading. For such a
purpose, it is enough to assume that the elements~$\sigma^{\nu}$ of the colour
group~$\cal G$ are superalgebra isomorphisms, or, in other words, that they
satisfy the condition
\begin{equation}
  \sigma^{\nu} \circ \gamma_q = \gamma_{q^{\nu}} \circ \sigma^{\nu},
  \label{eq:grading}
\end{equation}
in addition to Eq.~(\ref{eq:iso}). So we obtain
%
\begin{definition}   \label{def-colsuper}
A coloured Hopf algebra~${\cal H}^c$, as considered in
Definition~\ref{def-colHopf}, is called a coloured Hopf superalgebra if ${\cal H}_q$
is a Hopf superalgebra, and Eq.~(\ref{eq:grading}) is satisfied.
\end{definition}
\par
%
%
As a consequence of Eq.~(\ref{eq:permute}), if the Hopf superalgebras~${\cal H}_q$
are quasitriangular, then the corresponding universal $\cal R$-matrices ${\cal
R}_q$ satisfy the graded YBE~\cite{chaichian}. Hence, the coloured universal $\cal
R$-matrix ${\cal R}^c$, obtained from Definition~\ref{def-quasi}, is a solution of
the coloured graded~YBE.\par
%
%
Combining now Definitions~\ref{def-colQUEA} and~\ref{def-colsuper}, we may
consider coloured QUEAs of Lie superalgebras. As an example, let us consider the
two-parameter quantization~$U_{q,s}(gl(1/1))$ of the enveloping algebra
of $gl(1/1)$~\cite{burdik,dabrowski}.\par
%
%
The quantum superalgebra~$U_{q,s}(gl(1/1))$, for which $k = \C$, and $q$, $s \in \C
\setminus \{0\}$, is generated by two even generators $H$, $Z$, and two odd ones
$\psi^{\pm}$, with the relations~\cite{burdik}
\begin{equation}
  \left[H, \psi^{\pm}\right] = \pm 2 \psi^{\pm}, \qquad \left[Z, H\right] =
  \left[Z, \psi^{\pm}\right] = 0, \qquad  \left\{\psi^+, \psi^-\right\} =
  \frac{q^{2Z}-1}{q^2-1}, \qquad \left(\psi^{\pm}\right)^2 = 0,     \label{eq:comrel}  
\end{equation}
and a coalgebra structure depending upon both parameters $q$ and~$s$.\par
%
%
Eq.~(\ref{eq:comrel}) is left invariant under the grading-preserving transformations
\begin{equation}
  \sigma^{\nu}(H) = H, \qquad \sigma^{\nu}(Z) = \nu Z, \qquad
  \sigma^{\nu}(\psi^{\pm}) = a^{\nu} \psi^{\pm}, \qquad   
\end{equation}
where $\nu \in {\cal C} = \C \setminus \{0\}$ and $a^{\nu} \equiv \bigl((q^{2\nu}-1)
/ (q^2-1)\bigr)^{1/2}$, provided $q$ is changed into $q^{\nu}$ ($\nu$th power
of $q$), while $s$ is left unchanged. Hence, ${\cal Q} = \left(\C \setminus
\{0\}\right) \times \{s\}$, corresponding to a varying-parameter case. Since
$\nu' \circ \nu = \nu' \nu$, $\nu^0 = 1$, $\nu^i = \nu^{-1}$, the colour group~$\cal G$
is isomorphic to the abelian group
$Gl(1,\C)$.\par
%
%
The coloured maps and universal $\cal R$-matrix are easily 												
			obtained as
\begin{eqnarray}
  \Delta^{\lambda,\mu}_{q,s,\nu}\left(H\right) & = & H \otimes 1 + 1 \otimes H,
           \qquad \Delta^{\lambda,\mu}_{q,s,\nu}\left(Z\right) = \frac{\lambda}
           {\nu}\, Z \otimes 1 + \frac{\mu}{\nu}\, 1 \otimes Z, \nonumber \\
  \Delta^{\lambda,\mu}_{q,s,\nu}\left(\psi^{\pm}\right) & = & \frac{a^{\lambda}}
           {a^{\nu}} \psi^{\pm} \otimes s^{\mp\mu Z/2} q^{\mu Z} + \frac{a^{\mu}}
           {a^{\nu}} s^{\pm\lambda Z/2} \otimes \psi^{\pm}, \nonumber \\
  \epsilon_{q,s,\nu}(X) & = & 0, \qquad X \in \{H, Z, \psi^{\pm}\}, \nonumber \\
  S^{\mu}_{q,s,\nu}\left(H\right) & = & - H, \qquad
           S^{\mu}_{q,s,\nu}\left(Z\right) = - \frac{\mu}{\nu} Z, \qquad
           S^{\mu}_{q,s,\nu}\left(\psi^{\pm}\right) = - \frac{a^{\mu}}{a^{\nu}} q^{-\mu Z}
           \psi^{\pm}, \nonumber \\
  {\cal R}^{\lambda,\mu}_{q,s} & = & q^{(\mu H \otimes Z + \lambda Z \otimes H)/2}
           s^{(\mu H \otimes Z - \lambda Z \otimes H)/2} \nonumber \\
  & & \mbox{} \times \left\{1 \otimes 1
           - \sqrt{\left(q^{2\lambda}-1\right) \left(q^{2\mu}-1\right)}\, 
           s^{-\lambda Z/2} \psi^+ \otimes s^{-\mu Z/2} q^{-\mu Z} \psi^-\right\}. 
\end{eqnarray}
\par
%
%
The matrix representation of the coloured universal $\cal R$-matrix in any
finite-dimensional representation of $U_{q,s}(gl(1,1))$ provides us with a matrix
solution $R^{\lambda,\mu}_{q,s}$ of the coloured graded~YBE. For instance, in the
two-dimensional representation,
\begin{equation}
  D(Z) = \left(\begin{array}{cc}
             1 & 0 \\ 0 & 1
             \end{array} \right), \qquad  
  D(H) = \left(\begin{array}{cc}
             1 & 0 \\ 0 & -1
             \end{array} \right), \qquad 
  D(\psi^+) = \left(\begin{array}{cc}
                    0 & 1 \\ 0 & 0
                    \end{array} \right), \qquad
  D(\psi^-) = \left(\begin{array}{cc}
                    0 & 0 \\ 1 & 0
                    \end{array} \right),  
\end{equation}
we get the following $4 \times 4$ matrix solution of the coloured graded~YBE,
\begin{eqnarray}
  R^{\lambda,\mu}_{q,s} & \equiv & (D \otimes D) \left({\cal
        R}^{\lambda,\mu}_{q,s}\right) \nonumber \\
  & = &\left(\begin{array}{cccc}
        \ss q^{(\lambda+\mu)/2} s^{(-\lambda+\mu)/2} & \ss 0 & \ss 0 & \ss 0 \\
        \ss 0 & \ss q^{(-\lambda+\mu)/2} s^{(\lambda+\mu)/2} & \ss
             q^{-(\lambda+\mu)/2} \sqrt{\left(q^{2\lambda}-1\right)
             \left(q^{2\mu}-1\right)} & \ss 0 \\
        \ss 0 & \ss 0 & \ss q^{(\lambda-\mu)/2} s^{-(\lambda+\mu)/2} & \ss 0 \\
        \ss 0 & \ss 0 & \ss 0 & \ss q^{-(\lambda+\mu)/2} s^{(\lambda-\mu)/2}    
  \end{array}\right).   
\end{eqnarray}
The latter is related to the coloured $R$-matrix previously derived by
Burd\'\i k and Hellinger~\cite{burdik} by considering $2 \times 2$
representations of $U_{q,s}(gl(2))$ characterized by different eigenvalues
$\lambda$,~$\mu$ of~$Z$.\par
%
%

\end{document}